\documentclass[11pt,reqno]{amsart}
\usepackage{enumerate}
\usepackage{amsfonts}
\usepackage{amsmath}
\usepackage{amsthm}
\usepackage{amssymb}
\usepackage{latexsym}
\usepackage{multicol}
\usepackage{verbatim}
\usepackage{amscd}
\usepackage{hyperref}
\usepackage{pb-diagram}

\advance\textwidth by 1.2in \advance\oddsidemargin by -.6in
\advance\evensidemargin by -.6in
\newtheorem*{cor}{Corollary}
\newtheorem*{lem}{Lemma}
\newtheorem*{prop}{Proposition}

\theoremstyle{definition}
\newtheorem*{defn}{Definition}
\theoremstyle{definition}
\newtheorem{thm}{Theorem}

\newtheorem*{rem}{Remark}

\newenvironment{pf}{\proof}{\endproof}
\newcounter{cnt}
\newenvironment{enumerit}{\begin{list}{{\hfill\rm(\roman{cnt})\hfill}}{%
\settowidth{\labelwidth}{{\rm(iv)}}\leftmargin=\labelwidth%
\advance\leftmargin by
\labelsep\rightmargin=0pt\usecounter{cnt}}}{\end{list}}
\makeatletter
\def\mydggeometry{\makeatletter\dg@YGRID=1\dg@XGRID=20\unitlength=0.003pt\makeatother}
\makeatother \theoremstyle{remark}

\numberwithin{equation}{section}

\DeclareMathOperator{\Ht}{ht}

\let\bwdg\bigwedge
\def\bigwedge{{\textstyle\bwdg}}

\begin{document}

\newcommand{\thmref}[1]{Theorem~\ref{#1}}
\newcommand{\secref}[1]{Section~\ref{#1}}
\newcommand{\lemref}[1]{Lemma~\ref{#1}}
\newcommand{\propref}[1]{Proposition~\ref{#1}}
\newcommand{\corref}[1]{Corollary~\ref{#1}}
\newcommand{\remref}[1]{Remark~\ref{#1}}
\newcommand{\defref}[1]{Definition~\ref{#1}}
\newcommand{\er}[1]{(\ref{#1})}
\newcommand{\id}{\operatorname{id}}
\newcommand{\sgn}{\operatorname{sgn}}
\newcommand{\wt}{\operatorname{wt}}
\newcommand{\tensor}{\otimes}
\newcommand{\from}{\leftarrow}
\newcommand{\nc}{\newcommand}
\newcommand{\rnc}{\renewcommand}
\newcommand{\dist}{\operatorname{dist}}
\newcommand{\qbinom}[2]{\genfrac[]{0pt}0{#1}{#2}}
\nc{\cal}{\mathcal} \nc{\goth}{\mathfrak} \rnc{\bold}{\mathbf}
\renewcommand{\frak}{\mathfrak}
\newcommand{\supp}{\operatorname{supp}}
\renewcommand{\Bbb}{\mathbb}
\nc\bomega{{\mbox{\boldmath $\omega$}}} \nc\bpsi{{\mbox{\boldmath
$\Phi$}}}
 \nc\balpha{{\mbox{\boldmath $\alpha$}}}
 \nc\bpi{{\mbox{\boldmath $\pi$}}}

\newcommand{\lie}[1]{\mathfrak{#1}}
\makeatletter
\def\section{\def\@secnumfont{\mdseries}\@startsection{section}{1}%
  \z@{.7\linespacing\@plus\linespacing}{.5\linespacing}%
  {\normalfont\scshape\centering}}
\def\subsection{\def\@secnumfont{\bfseries}\@startsection{subsection}{2}%
  {\parindent}{.5\linespacing\@plus.7\linespacing}{-.5em}%
  {\normalfont\bfseries}}
\makeatother
\def\subl#1{\subsection{}\label{#1}}
 \nc{\Hom}{\operatorname{Hom}}
\nc{\End}{\operatorname{End}} \nc{\wh}[1]{\widehat{#1}}
\nc{\Ext}{\operatorname{Ext}} \nc{\ch}{\text{ch}}
\nc{\ev}{\operatorname{ev}} \nc{\Ob}{\operatorname{Ob}}
\nc{\soc}{\operatorname{soc}} \nc{\rad}{\operatorname{rad}}
\nc{\head}{\operatorname{head}}
\def\Im{\operatorname{Im}}
\def\gr{\operatorname{gr}}
\def\mult{\operatorname{mult}}

 \nc{\Cal}{\cal} \nc{\Xp}[1]{X^+(#1)} \nc{\Xm}[1]{X^-(#1)}
\nc{\on}{\operatorname} \nc{\Z}{{\bold Z}} \nc{\J}{{\cal J}}
\nc{\C}{{\bold C}} \nc{\Q}{{\bold Q}}
\renewcommand{\P}{{\cal P}}
\nc{\N}{{\Bbb N}} \nc\boa{\bold a} \nc\bob{\bold b} \nc\boc{\bold c}
\nc\bod{\bold d} \nc\boe{\bold e} \nc\bof{\bold f} \nc\bog{\bold g}
\nc\boh{\bold h} \nc\boi{\bold i} \nc\boj{\bold j} \nc\bok{\bold k}
\nc\bol{\bold l} \nc\bom{\bold m} \nc\bon{\bold n} \nc\boo{\bold o}
\nc\bop{\bold p} \nc\boq{\bold q} \nc\bor{\bold r} \nc\bos{\bold s}
\nc\bou{\bold u} \nc\bov{\bold v} \nc\bow{\bold w} \nc\boz{\bold z}
\nc\boy{\bold y} \nc\ba{\bold A} \nc\bb{\bold B} \nc\bc{\bold C}
\nc\bd{\bold D} \nc\be{\bold E} \nc\bg{\bold G} \nc\bh{\bold H}
\nc\bi{\bold I} \nc\bj{\bold J} \nc\bk{\bold K} \nc\bl{\bold L}
\nc\bm{\bold M} \nc\bn{\bold N} \nc\bo{\bold O} \nc\bp{\bold P}
\nc\bq{\bold Q} \nc\br{\bold R} \nc\bs{\bold S} \nc\bt{\bold T}
\nc\bu{\bold U} \nc\bv{\bold V} \nc\bw{\bold W} \nc\bz{\bold Z}
\nc\bx{\bold x}
\title[]{ Ideals in parabolic subalgebras of simple Lie algebras}
\author{Vyjayanthi Chari, R. J. Dolbin and T. Ridenour}
\thanks{VC was partially supported by the NSF grant DMS-0500751}
\address{Department of Mathematics, University of
California, Riverside, CA 92521.} \email{chari@math.ucr.edu}
 \email{rjdolbin@math.ucr.edu}\email {tbr4@math.ucr.edu}
\maketitle
\begin{abstract}{We study ad--nilpotent ideals of a parabolic subalgebra of a simple Lie algebra. Any such ideal determines an antichain in a set of positive roots of the simple Lie algebra. We give a necessary and sufficient condition for an antichain to determine an ad--nilpotent ideal of the parabolic. We write down all such antichains for the classical simple Lie algebras and in particular recover the results of D. Peterson. In section 2 of the paper we study the unique ideal in a parabolic which is irreducible as a module for the reductive part and give several equivalent statements that are satisfied by the corresponding subset of roots.   }\end{abstract}

\section*{Introduction} In recent years, there have been a number of articles  \cite{CP1}, \cite{CP2}, \cite{Ko}, \cite{Pan}, \cite{Sut} (to name just a few)  on ad--nilpotent ideals in a Borel subalgebra of a simple Lie algebra. These papers were motivated by a result of D. Peterson who showed that there are exactly $2^n$ such ideals where $n$ is the rank of $\lie g$. In \cite{CP1} and \cite{CP2}, a bijection is established between the set of all ideals and a certain subset of elements of the affine Weyl group. This method was later generalized in \cite{Rig}, \cite{Rig2}  and used to study ad--nilpotent ideals in a parabolic subalgebra $\lie p$ of $\lie g$.

In section one of this paper, we approach this problem from a more elementary perspective. Let $\lie b$ be a Borel subalgebra of a simple Lie algebra  $\lie g$ and  $R^+$ the corresponding set of positive roots.   Assume that $R^+$  is partially  ordered as usual: $\alpha\le \beta$ iff $\beta-\alpha$ is in the non--negative integers span of $R^+$. An antichain in $R^+$ is a subset of elements which are pairwise unrelated in this order. Assume now that $\lie p$ is a parabolic subalgebra of $\lie g$ containing $\lie b$.
Any ad--nilpotent ideal $\lie i$ in $\lie p$  is contained in the unipotent radical of $\lie p$  and hence determines a subset of $R^+$, and we let $A(\lie i)$ be the  antichain consisting of the minimal elements of this subset. We give a necessary and sufficient condition for an antichain to determine an ideal in the unipotent radical of $\lie p$. As an application of this condition, we write down explicitly  all the antichains which determine abelian ad--nilpotent ideals in a parabolic subalgebra of a simple Lie algebra of classical type. Our methods work efficiently in the exceptional cases as well although we do not list the antichains  for the exceptional algebras in this paper. In particular, we can count the number of ad--nilpotent  ideals in $\lie p$, this recovers the result of D. Peterson when $\lie p$ is just the Borel subalgebra and the results of \cite{Rig}, \cite{Rig2} for a general parabolic subalgebra.

In section two of this paper, we focus our attention on a particular family of subsets of $R$.
Any parabolic subalgebra $\lie p$ contains a unique ad--nilpotent  ideal $\lie i_0$ which is irreducible as a module for the Levi factor of  the parabolic.  Let $\Phi\subset R$ be the corresponding subset of roots. We give several different characterizations of such sets. We show that if $2\rho_\Phi=\sum_{\alpha\in\Phi}\alpha$  then $\Phi=\{\alpha\in R: (\rho_\phi,\alpha)\ {\rm{is\  maximal}}\}$. 
Equivalently, we prove that  if $\Phi\subset R$ is such that $(\alpha,\beta)\ge 0$ and for all  $\gamma,\delta\in R$ \ \  $\gamma+\delta\in\Phi+\Phi$ only if  $\gamma,\delta\in\Phi$, then there exists a parabolic subalgebra $\lie p$ of $\lie g$ such that  the ideal $\lie i_0$ is given by $\Phi$.
  As a consequence we can write down the subsets $\Phi(\lie i_0)$ explicitly. Section 2 is motivated by the results of \cite{CG}, where we prove that to each such subset one can define an infinite--dimensional associative algebra which is Koszul and of global dimension equal to the cardinality of the subset.
\subsection*{Acknowledgments}
We are grateful to Jacob Greenstein for many discussions. The first author also thanks  Bertram Kostant for pointing out the connection between the main theorem in Section 2 and his result on the decomposition of the unipotent part of a parabolic.

\section{Ad-nilpotent abelian ideals of a parabolic subalgebra}{}

\subsection{} Throughout the paper $\bz$ denotes the set of integers and $\bz_+$ the set of non--negative integers. Let $\lie g$ be a finite--dimensional complex simple Lie algebra
of rank $n$. Fix   a  Cartan subalgebra $\lie h$ of $\lie g$ and  let $R\subset\lie h^*$ be the corresponding root system and $W$ the Weyl group. Given $\alpha\in R$, let $\lie g_\alpha$ be the corresponding root space and fix elements $x_\alpha\in\lie g_\alpha$ so that $\lie g_\alpha=\bc x_\alpha$. Recall that if $\alpha,\beta\in R$ is such that $\alpha+\beta\in R$ then $[x_\alpha,x_\beta]=c x_{\alpha+\beta}$ for some non--zero $c\in \bc$.

Let $(\ ,\ )$ be the symmetric non--degenerate form on $\lie h^*$ which is induced  from  the restriction to $\lie h$ of the Killing form of  $\lie g$.
Set $I=\{1,\cdots,n\}$ and let $\{\alpha_i:i\in I\}$ (resp. $\{\omega_i: i\in I\}$)
  be a set of simple roots (resp. fundamental weights) and
  $Q$ and  $P$   (resp. $Q^+$, $P^+$) be the $\bz$--span (resp. $\bz_+$--span) of the simple roots and  fundamental weights  respectively.  Set $R^+=R\cap Q^+$

For $i\in I$, define
      $d_i: Q\to\bz$ by requiring:$$\eta=\sum_{i\in I}d_i(\eta)\alpha_i,$$
and set $\Ht\eta=\sum_{i\in I} d_i(\eta)$. Clearly $d_i$ and $\Ht$ are additive homomorphisms of abelian groups.

Define a partial order on $P$ by $\lambda\le \mu$ iff $\mu-\lambda\in Q^+$. If $\theta$ is the highest root of $R^+$, then $\theta$
 is the unique maximal element in $R^+$ with respect to this
 order. Given $\lambda,\mu\in P$ with $\lambda\le \mu$, let
 $k(\mu, \lambda)$ be the minimal non--negative integer so that
 $\mu-\lambda$ can be written as a sum of $k(\lambda,\mu)$ (not
 necessarily distinct)
 elements of $R^+$.
We shall need the following elementary lemma.
\begin{lem}\label{lemp}
\begin{enumerit}\item[(i)]
 Let $\lambda\in P^+$ and assume that $\lambda\le\theta$.
 Then either $\lambda=0$ or $\lambda\in R$.
 \item[(ii)] Suppose that $\alpha,\beta\in R^+$ and that $\beta<\alpha$. There exists $\gamma,\delta\in R^+ $ with $\gamma\ge\beta$ such that $\alpha=\gamma+\delta$.
     \item[(iii)] Suppose that $\alpha,\beta,\gamma\in R$ are such that $\alpha+\beta$, $\alpha+\beta+\gamma\in R$. Then either $\alpha+\gamma\in R$ or $\beta+\gamma\in R$.
         \end{enumerit}
 \end{lem}
 \begin{pf} To prove (i) note that $\theta$ is the highest weight of the adjoint representation of $\lie g$.
  It follows from \cite{Hum} that if $\lambda\le\theta$ then $\lambda$ must be a weight of the adjoint representation of $\lie g$ which proves (i).

    We proceed by induction on $k(\alpha,\beta)$ with induction obviously beginning when $k(\alpha,\beta)=1$.
    Assume we have proved the result for all pairs $\gamma$,
$\gamma'$ with $\gamma>\gamma'$ and $k(\gamma,\gamma')<r$.
  Let $\alpha>\beta$ be such that $k(\alpha,\beta)=r$ and write $\alpha=\beta+\beta_1+\cdots +\beta_r$ for some $\beta_p\in R^+$, $1\le p\le r$.
   The minimality of $k(\alpha,\beta)$ implies immediately that $\beta_s+\beta_p\notin R^+$
    and hence $(\beta_s,\beta_p)\ge 0$ for all $1\le s,p\le r$.
     If $\beta+\beta_r\notin R^+$ then we have $(\beta,\beta_r)\ge 0$
     and hence $(\alpha,\beta_r)\ge (\beta_r,\beta_r)>0$.
  Since $\alpha\ne\beta_r$, it follows that    $\alpha'=\alpha-\beta_r\in R^+$. The result follows since $\alpha'\ge
  \beta$.    If $\beta+\beta_r\in R^+$,
        then   $k(\alpha,\beta+\beta_r)<r$.
         Hence the inductive hypothesis gives $\alpha=\gamma+\delta$ for some $\gamma,\delta\in R^+$ with $\gamma \ge \beta + \beta_r > \beta$ and part (ii) is proved.

To prove (iii),
 note that $x_{\alpha+\beta+\gamma}= c[x_\gamma, x_{\alpha+\beta}]= d[x_\gamma,[x_\alpha,x_\beta]]$
  for some non-zero $c,d\in\bc$. The result is immediate from the Jacobi identity.
 \end{pf}

\subsection{} Given $J\subset I$ set $$R(J)=\{\alpha\in R: d_i(\alpha)= 0 \ {\rm{if}}\ i\notin J\},\ \ R^+(J)=R(J)\cap R^+.$$
 \begin{defn} A subset $\Phi$ of $ R^+$
  is called a $J$--ideal if  \ $\Phi\cap R^+(J)=\emptyset$ and
$$ \alpha\in\Phi,\ \beta\in R^+\cup R(J),\ \beta+\alpha\in R\implies \beta+\alpha\in \Phi.$$

 A subset $A$ of $R^+$ is called a $J$--antichain if $A\cap R^+(J)=\emptyset$ and for all $\alpha,\beta\in A$ and $j\in J$, we have $ \alpha\nleq\beta$, $\beta\nleq\alpha$ and
$\alpha-\alpha_j\notin R$.  \hfill\qedsymbol\end{defn}
If  $J\subset J'$  then any $J'$--ideal (resp. $J'$--antichain) is also a $J$--ideal (resp. $J$--antichain). In the case when $J=\emptyset$  we drop the dependence on $J$, for instance an ideal is  a $\emptyset$-ideal.
 \begin{rem}  Let $\lie g_\alpha$ be the root space of $\lie g$ corresponding to $\alpha$. Then $\Phi$ is a $J$--ideal iff  the subspace $\oplus_{\alpha\in \Phi}\lie g_\alpha$ is an ad--nilpotent ideal of the parabolic subalgebra $\lie p_J=\lie h\oplus_{\alpha\in R^+}\lie g_\alpha \oplus_{\alpha\in R^+(J)}\lie g_{-\alpha}$ of $\lie g$. Conversely  any ad-nilpotent ideal in $\lie p_J$ determines a $J$--ideal in $R^+$.  The set of $J$--antichains determine a minimal set of generators of the ideal in $\lie p_J$.\end{rem}

\subsection{} It is immediate from the definition that, if $A$ is a $J$--antichain, then\begin{equation}\label{inn}\alpha,\beta\in A, \alpha\ne\beta\implies (\alpha, \beta)\le 0,\ \ \alpha\in A, j\in J\implies (\alpha,\alpha_j)\le 0.\end{equation}
 \begin{lem} \label{equiv} Let $A$ be a $J$--antichain in $R^+$. For all $\alpha\in A$, $\gamma\in R^+(J)$ we have $\alpha-\gamma\notin R$.
 \end{lem}
 \begin{pf} The Lemma is proved by an induction on $\Ht \gamma$ with induction beginning at  $\Ht(\gamma)=1$ since $A$ is $J$--antichain. Suppose that we have proved the Lemma for all $\gamma'\in R^+(J)$ with $\Ht(\gamma')<r$. Let $\gamma\in R^+(J)$ be such that $\Ht(\gamma)=r$ and choose $j\in J$ such that $(\gamma,\alpha_j)>0$. Then $\gamma-\alpha_j\in R^+(J)$.  Suppose that  $\alpha\in A$ is such that $\alpha-\gamma\in R$. Then $(\alpha-\gamma,\alpha_j)<0$ by \eqref{inn} and, hence, we have $\alpha-\gamma+\alpha_j\in R$. But this contradicts the induction hypothesis since $\Ht(\gamma-\alpha_j)=r-1$.
 \end{pf}

 \subsection{} Let  $\Phi$ be a $J$--ideal and let $A(\Phi)$ be the set of minimal elements of $\Phi$ with respect to the partial order $\le$. Conversely, if $A$ is $J$--antichain, set $$\Phi(A)=\bigcup_{\beta\in A}\{\alpha\in R^+: \alpha\ge\beta\ \}.$$
\begin{prop} The assignment $A\to\Phi(A)$ is a  bijection between the set of  $J$--antichains in  $R^+$ and  $J$--ideals in $R^+$.
\end{prop}
\begin{pf}

 Suppose that $\beta\in A$  and that $\alpha\ge \beta$; then, $d_i(\alpha)\ge d_i(\beta)$ for all $i\in I$.  Since $\beta\notin R^+(J)$ there exists $i_0\in I\setminus J$ such that $d_{i_0}(\alpha)\ge d_{i_0}(\beta)>0$ and  hence $\alpha\notin R^+(J)$. If $\gamma\in R^+$ is such that $\alpha+\gamma\in R^+$ then clearly $\alpha+\gamma\ge\beta $ and hence $\alpha+\gamma\in \Phi$. Finally, suppose that $\gamma\in R^+(J)$ is such that $\alpha-\gamma\in R$.  We proceed by induction on $k(\alpha,\beta)$.  If $k(\alpha,\beta)=0$, then $\alpha=\beta$ and, by Lemma \ref{equiv}, we see that $\beta-\gamma\notin R$ and there is nothing to prove.  Assume that we have proved the result for all $\alpha'\in R^+$ with $\alpha'\ge\beta$ and $k(\alpha',\beta)< r$. If $\alpha\ge\beta$ and $k(\alpha,\beta)=r$, using Lemma \ref{lemp}(ii) we can write $\alpha=\alpha'+\gamma'$ where $\alpha'\ge\beta$ and $\gamma'\in R^+$ and hence $k(\alpha',\beta)<k(\alpha,\beta)$.
 Applying Lemma \ref{lemp}(iii)  to $\alpha'+\gamma'-\gamma$ it follows that either $\alpha'-\gamma\in R$ or $\gamma'-\gamma\in R$. In the first case the induction hypothesis applied to $\alpha'$ gives  $\alpha>\alpha'-\gamma\ge\beta$ and we  are done. If $\gamma'-\gamma\in R^+$ we are done since again we have $\alpha-\gamma\ge\alpha'\ge\beta$. If $\gamma'-\gamma\in R^-$ then in fact $\gamma\in R^+(J)$ and we have by the induction hypothesis again that  $\alpha-\gamma=\alpha'-(\gamma-\gamma')\ge \beta$ and the Proposition is proved.
\end{pf}

\subsection{} An ideal $\Phi$ in  $R^+$ is said to be of nilpotence $k$ if given  elements $\beta_1,\cdots,\beta_{k+1}\in \Phi$, (not necessarily
distinct), we have $\sum_{p=1}^{k+1}\beta_p\notin R$. The following result gives a necessary and sufficient condition for an  antichain to determine an ideal of nilpotence $k$.
\begin{thm}\label{antichain}  Let $A$ be  an antichain in $R^+$. The ideal $\Phi
(A)$ is  of nilpotence $k$ iff given $\beta_1,\cdots,\beta_{k+1}\in A$ (not necessarily distinct), we have $\sum_{s=1}^{k+1}\beta_s\nleq\theta$.
\end{thm}
\begin{pf} Suppose that $\Phi(A)$ is an ideal of nilpotence $k$ and assume for a contradiction that there exist $\beta_1,\cdots,\beta_{k+1}\in A$ such that $\sum_{s=1}^{k+1}\beta_s\leq\theta$. We claim that there exist elements $\gamma_j\in R^+$ with $\gamma_j\ge\beta_j$ for $1\le j\le k+1$ such that $\sum_{j=1}^{k+1}\gamma_j\in R^+$. Since  $\gamma_j\in\Phi(A)$ for all $1\le j\le k+1$ the claim implies that $\Phi(A)$ is not of nilpotence $k$ which is a contradiction.

 To see that induction begins, suppose that  $\sum_{s=1}^{k+1}\Ht\beta_s=\Ht\theta$. Since $\theta-\sum_{s=1}^{k+1}\beta_s\in Q^+$,  it follows that $\sum_{s=1}^{k+1}\beta_s=\theta$ and the claim follows by taking $\gamma_j=\beta_j$.
For the inductive step assume that we have proved the proposition for all $\beta_j$, $1\le j\le k+1$ with $\Ht\sum_{s=1}^{k+1}\beta_s >r$. If $\sum_{j=1}^{k+1}\beta_j\in R^+$, there is nothing to prove. If $\beta=\sum_{j=1}^{k+1}\beta_j\notin R^+$ then using Lemma \ref{lemp}(i), we see that the condition that $\beta\le\theta$  implies that $\beta\notin P^+$.
Since $\beta\ne 0$,   there exists $i_0\in I$  with $(\beta,\alpha_{i_0})<0$, in particular
 for some $1\le s\le k+1$ we must have $(\beta_s,\alpha_{i_0})<0$, i.e $\beta_s+\alpha_{i_0}\in R^+$. The inductive step follows if we prove that $\sum_{j=1}^{k+1}\beta_j+\alpha_{i_0}\le\theta$. This is clear by noting that since $$(\theta- \sum_{j=1}^{k+1}\beta_j, \alpha_{i_0})>0,$$
we must have $$\theta-\sum_{j=1}^{k+1}\beta_j=\sum_{i\in I}r_i\alpha_i,\ \ r_{i_0}\ne 0.$$ The converse direction of the theorem is clear.
\end{pf}

\subsection{} An ideal of nilpotence one is called an abelian ideal and an abelian antichain is one that defines an abelian ideal. If $J$ is any subset of $I$ we have analogous notions of abelian  $J$--ideals and  abelian $J$--antichains.

The following is an immediate consequence of Theorem \ref{antichain}.
\begin{prop} \label{abelian} Let $J$ be any subset of $I$.  Let $A$ be a $J$--antichain  in   $R^+$ . Then $A$ is an abelian  $J$--antichain iff the following holds:
\begin{enumerit}
\item Given $\alpha,\beta\in A$ with $\alpha\ne\beta$,   there exists $i\in I$ (depending on $\alpha$, $\beta$) such that $d_i(\alpha)+d_i(\beta)>d_i(\theta)$, in particular $d_i(\alpha)\ne 0$ and $d_i(\beta)\ne 0$.
\item Given  $\alpha\in A$, there exists $i\in I$ such that $2d_i(\alpha)>d_i(\theta)$.\end{enumerit}
\hfill\qedsymbol
\end{prop}

Let $J\subset I$. For $s\ge 1$, let $\ba_{s, J}$ be the set of  $J$--abelian antichains with $s$ elements and set $\ba_{0,J}=\emptyset$. As an application of the proposition we write down all elements of $\ba_{s,J}$, $s\ge 1$ for the classical Lie algebras. We shall assume that the set of simple roots of $\lie g$ is numbered as in \cite{B}.
 We compute the cardinality of $\ba_{s,J}$ for each $s\ge 0$. In the case when $J=\emptyset$, we    prove that $\sum_{k\ge 0}\# \ba_k=2^n$  thus recovering a theorem D. Peterson, \cite{Ko},\cite{Sut}. The case when $J\ne\emptyset$ recovers  the results of \cite{Rig}.
 From now on we set $\ba_s=\ba_{s,\emptyset}$.

\subsection{} Henceforth, we shall understand that $\binom{n}{k}=0$ if $k>n$. Let  $R$ be of type $A_n$. For $i,j\in I$,  set $\alpha_{i,j}=\alpha_i+\alpha_{i+1}+\cdots+\alpha_j$  and note that $\alpha_{i,i}=\alpha_i$. Then $$R^+=\{\alpha_{i,j}:i,j\in I, i\le j\},\ \ \theta=\alpha_{1,n}.$$
\begin{prop} For $s\ge 1$, we have \begin{equation}\label{ans}\ba_s=\{\{\alpha_{i_k,j_k}\}_{ 1\le k\le s}: \ i_k,j_k\in I, i_k<i_{k+1}, j_k<j_{k+1}, \ i_s\le j_1\}
\end{equation}
In particular, $$\#\ba_s=\binom{n}{2s}+\binom{n}{2s-1},\ \ \ \  \sum_{s\ge 0} \#\ba_s=2^n.$$
\end{prop}
\begin{pf} If $s=1$ the statement is immediate from Proposition \ref{abelian} since $2\alpha\nleq\theta$ for all $\alpha\in R^+$.
Set
 $A=\{\alpha_{i_k,j_k}:  1\le i_1<i_2<\cdots <i_s\le j_1<j_2<\cdots <j_s\le n\}$. Since \begin{gather*}d_{i_r}(\alpha_{i_r,j_r}-\alpha_{i_{r+1},j_{r+1}})=\ 1\qquad
  d_{j_{r+1}}(\alpha_{i_r,j_r}-\alpha_{i_{r+1},j_{r+1}})=-1\end{gather*} it follows that  $A$ is an antichain.  Moreover, for all $1\le r\le q\le s$, we have $$d_{i_{q}}(\alpha_{i_r,j_r}+\alpha_{i_{q},j_{q}})= 2>d_{i_q}(\theta),\ \ $$ and hence Proposition \ref{abelian} implies that $A\in\ba_s$. Conversely, suppose that $A\in\ba_s$ for some $s>1$, say $A=\{\alpha_{i_1,j_1},\cdots,\alpha_{i_s,j_s}\}$. If   $i_r=i_{r+1}$ (resp. $j_r=j_{r+1}$) for some $1\le r\le s-1$, then  $\alpha_{i_r,j_r}-\alpha_{i_r,j_{r+1}}\in R$  (resp. $\alpha_{i_r,j_r}-\alpha_{i_{r+1},j_r}\in R$) contradicting the fact that  $A$ is an antichain. Hence we can assume without loss of generality that $i_1<i_2< \cdots < i_s$. If $j_r<j_{r-1}$  for some $r\ge 2$, then we have $\alpha_{i_r,j_r}<\alpha_{i_{r-1},j_{r-1}}$ which again contradicts the fact that $A$ is an antichain. Hence, we must have $j_1<\cdots< j_s$. Finally note that if $i_s>j_1$ then $\alpha_{i_1,j_1}+\alpha_{i_s,j_s}\le\theta$ which is impossible by Proposition \ref{abelian} and hence \eqref{ans} is proved. The final statement of the proposition is clear with  the two terms in the $\#\ba_s$ coming from the case when $i_s<j_1$ and $i_s=j_1$ respectively.
\end{pf}
The following is now immediate from the definition of a $J$--antichain together with the remark that $\alpha_{k,\ell}-\alpha_{k,k}\in R$ for all $\ell\ne k$..
\begin{cor} Suppose that $J$ is a  subset of $I$. Then $$\ba_{s,J}=\{\{\alpha_{i_k,j_k}\}_{ 1\le k\le s} \in \ba_s\ :\ i_k,j_k\notin J\}$$ 
In particular, $$\#\ba_{s,J}=\binom{n-\#J}{2s}+\binom{n-\#J}{2s-1},\ \ \ \  \sum_{s\ge 0} \#\ba_{s,J}=2^{n-\#J}.$$
\end{cor}

\subsection{} If $R$ is of type $B_n$  set
\begin{gather*}\alpha_{i,j}=\alpha_i+\cdots+\alpha_j,\ \ i,j\in I,\ \ \ \  \beta_{k,\ell}=\alpha_{k,n}+\alpha_{\ell,n},\ \ k,\ell\in I, k\ne\ell,\end{gather*}
and then,  $ R^+=\{\alpha_{i,j}: i, j\in I, i\le j\}\cup\{ \beta_{k,\ell} : \  k,\ell\in I, k<\ell\}, \ \ \theta=\beta_{1,2}.$
\begin{prop} For $s\ge 1$, we have  $\ba_s=\ba_s^1\sqcup\ba_s^2$ where \begin{gather*}
\ba_s^1=\{\{\beta_{i_k,j_k}\}_{1\le k\le s}:\ \  i_k,j_k\in I,\ \  i_k<i_{k+1},\ \ j_{k+1}<j_k, \ \  i_s<j_s\}\\
 \ba_s^2= \{\{\alpha_{1,\ell},\{\beta_{i_k,j_k}\}_{:1\le k\le s-1}: \ \ \{\beta_{i_k,j_k}\}_{:1\le k\le s-1}\in\ba_{s-1}^1,\ \ell\in I,\ i_1\ge 1 \  j_1\le \ell\}
  \end{gather*}
In particular $$\#\ba^1_s= \binom{n}{2s}, \qquad \#\ba_s^2= \binom{n-1}{2s-2}+\binom{n-1}{2s-1},\  \qquad \sum_{s \ge 0} \#\ba_s=2^n.$$
\end{prop}
\begin{pf} Suppose that $1\le i_p<i_q<j_q<j_p\le n$.  The equations \begin{gather*} d_{i_p}(\beta_{i_p,j_p}-\beta_{i_q,j_q})=1,\ \ d_{j_q}(\beta_{i_p,j_p}-\beta_{i_q,j_q})=-1\\
d_{n}(\beta_{i_p,j_p}+\beta_{i_q,j_q})=4= 2d_n(\beta_{i_p,j_p})\end{gather*} prove along with Proposition \ref{abelian} that $\ba_s^1\subset \ba_s$. If $i_p>1$ and $j_p<\ell$, then the fact that $\ba_s^2\subset \ba_s$  follows by using the  additional equations
\begin{gather*} d_1(\alpha_{1,\ell}-\beta_{i_p,j_p})=1, \ \ d_n(\alpha_{1,\ell}-\beta_{i_p,j_p})\le -1,\\
2d_1(\alpha_{1,\ell})=2,\ \ d_{j_p}(\alpha_{1,\ell}+\beta_{i_p,j_p})=3.\end{gather*}

 For the converse let $A\in\ba_s$.  Suppose first that $A=\{\beta_{i_1,j_1},\cdots ,\beta_{i_s,j_s}\}$. As in the case of $A_n$, the fact that $A$ is an antichain means we can assume without loss of generality that $i_1<i_2<\cdots < i_s$ which in turn forces $j_s<\cdots <j_2<j_1$. This proves that $A\in\ba^1_s$.

 Suppose now that $\alpha_{i_0,j_0}\in A$ for some $1\le i_0\le j_0\le n$. Since $$2\alpha_{1,q}\nleq\theta,\ \  \,\ \ 2\alpha_{p,q}\le \theta,\ \  p>1,$$
it follows from   Proposition \ref{abelian} that  $i_0=1$. Since the elements
   $\alpha_{1,q}$ and $\alpha_{1,p}$ are always related it follows that  there exists a unique $\ell\in I$ such that $\alpha_{1,\ell}\in A$ and hence  the set $A\setminus \{\alpha_{1,}\}\in\ba_{s-1}^2$, and so we have, $$A\setminus\{\alpha_{1,\ell}\} =\{\beta_{i_1,j_1},\cdots ,\beta_{i_{s-1},j_{s-1}}\}, \ \  1 \le i_1<i_2<\cdots i_{s-1}<j_{s-1}<\cdots <j_1\le n.$$ Since $\alpha_{1,\ell}\le \beta_{1, j}$ for all $j\ge 2$ it follows that we must have $i_1>1$. To prove that $j_1<\ell$ it suffices to note that $\alpha_{1,j_0}+\beta_{i_1,j_1}\le\theta$ if $j_1<j_0$. This completes the proof of the proposition.

\end{pf}
\begin{cor} Suppose that $J$ is a subset of $I$. Then $\ba_{s,J}=\ba^{1,1}_{s,J}\sqcup\ba^{1,2}_{s,J}\sqcup \ba^2_{s,J}$ where
 \begin{gather*}\ba^{1,1}_{s,J}=\{\{\beta_{i_k,j_k}\}_{1\le k\le s}\in\ba_s^1: i_k, j_k\notin J \},\\
 \ba^{1,2}_{s,J}=\{\{\beta_{i_k,j_k}\}_{1\le k\le s}\in\ba_s^1: \{\beta_{i_k,j_k}\}_{1\le k\le s-1}\in\ba^{1,1}_{s-1,J},\ \  i_s\in J, j_s=i_s+1\notin J\ \},\\
 \ba_{s,J}^{2,1}=
 \{\{\alpha_{1,\ell},\{\beta_{i_k,j_k}\}_{1\le k\le s-1}\}: \ell\notin J,\ \ \{\beta_{i_k,j_k}\}_{1\le k\le s-1}\in\ba^1_{s-1,J}\},\qquad \ 1\notin J,
  \end{gather*}
  and $\ba_{s,J}^2=\emptyset$ if $1\in J$.

\end{cor}

\subsection{} If $R$ is of type $C_{n}$, set
\begin{gather*}
\alpha _{i,j}=\alpha _{i}+\cdots +\alpha _{j},\ \ i,j\in I,\ \
\beta _{k,\ell }=\alpha _{k,n-1}+\alpha _{\ell ,n},\ \ k,\ell\in I\setminus\{n\}.
\end{gather*}
and then, $
R^{+}=\{\alpha _{i,j}: \ \ i,j\in I,\ \ i\le j\}\  \cup\ \{ \beta _{k,\ell }: k,\ell\in I\setminus\{n\},\ \ k\le \ell\},\ \ \theta =\beta _{1,1}.
$ An analysis entirely analogous to the one for $B_n$ gives the following proposition and we omit the details.
\begin{prop} For  $s\ge 1$, we have  $\ba_s=\ba_s^1\sqcup\ba_s^2$ where
 \begin{gather*}\ba_s^1=\{
\{\beta _{i_{k},j_{k}}\}_{1\le k\le s}: \ i_k, j_k\in I\setminus\{n\}, ,\ \ i_k<i_{k+1},\ \ j_{k+1}<j_k,\ \ \ i_s\le j_s \},\\
\ba_s^2=\{\{\alpha_{\ell,n},\{\beta _{i_{k},j_{k}}\}_{1\le k\le s-1}: \{\beta _{i_{k},j_{k}}\}_{1\le k\le s-1}\in\ba_{s-1}^1, \ \  \ell\in I,\ \ \ell<i_1\}.\end{gather*}
In particular, $$\#\ba^1_s = \binom{n-1}{2s}+\binom{n-1}{2s-1},\qquad \#\ba_s^2=\binom{n-1}{2s-1}+\binom{n-1}{2s-2},\qquad \sum_{s\ge 0} \#\ba_s=2^n.$$\hfill\qedsymbol
\end{prop}
\begin{cor} Suppose that $J$ is a  subset of $I$. Then $\ba_{s,J}=\ba_{s,J}^1\sqcup\ba_{s,J}^2$, where \begin{gather*}\ba_{s,J}^1=\{\{\beta _{i_{k},j_{k}}\}_{1\le k\le s}\in\ba_s^1: i_k,j_k\notin J\},\ \ \\ \ba_{s,J}^2=
\{\alpha_{\ell,n},\{\beta _{i_{k},j_{k}}\}_{1\le k\le s-1}\}\in\ba_s^2: \ell\notin J,\ \{\beta _{i_{k},j_{k}}\}_{1\le k\le s-1}\in\ba_{s-1,J}^1\}, \  n\notin J\end{gather*} and $\ba_{s,J}^2=\emptyset$ otherwise. Hence $$\sum_{s\ge 0}\#\ba_{s,J}=2^{n-\#J}.$$
\hfill\qedsymbol

\end{cor}

\subsection{} The case of $D_n$ is analyzed in a similar way. It is, however, more tedious and the nice patterns for the cardinality of the sets $\ba_s$ are
broken. Set $\tilde I= I\setminus\{n-1,n\}$. The elements of $R^+$ are:
\begin{gather*}\alpha_{i,j}=\alpha_i+\cdots+\alpha_j,\qquad i,j\in\tilde I,\ \ i\le j,\\
\beta_{p,q}=\alpha_{p,n-2}+\alpha_{n-1}+\alpha_n+\alpha_{q,n-2},\qquad
p,q\in \tilde I\},\ p<q,\\
\gamma_{i,n-1}=\alpha_{i,n-2}+\alpha_{n-1},\ \
\gamma_{i,n}=\alpha_{i,n-2}+\alpha_n, \ \
\delta_i=\alpha_{i,n-2}+\alpha_{n-1}+\alpha_n,\ ,i\in \tilde I,
\end{gather*} and  $\theta=\beta_{1,2}$.

\begin{prop}  For $s\ge 1$ we have $\ba_s=\bigsqcup_{p=1}^6\ba_{s,J}^p$ where
   \begin{gather*} \ba_s^1=\{
\{\beta _{i_{k},j_{k}}\}_{1\le k\le s}: \ i_k, j_k\in \tilde I, ,\ \ i_k<i_{k+1},\ \ j_{k+1}<j_k,\ \ \ i_s< j_s \},\\
\ba_s^2=
\{\{\alpha_{1, j_0}, \{\beta _{i_{k},j_{k}}\}_{1\le k\le s-1}\}: \{\beta _{i_{k},j_{k}}\}_{1\le k\le s-1}\}\in\ba_{s-1}^1,\ j_0\in \tilde I, j_1<j_0
 \},\\
\ba_s^3=\{\{\gamma_{1,n-1},\gamma_{1,n},\{\beta _{i_{k},j_{k}}\}_{1\le k\le s-2}\}\}: \{\beta _{i_{k},j_{k}}\}_{1\le k\le s-2}\}\in\ba_{s-2}^1,\ \},\\\ba_s^4=\{\{\gamma_{i_0},\delta_{i_1}, \{\beta _{i_{k},j_{k}}\}_{2\le k\le s-1}\}: \{\beta _{i_{k},j_{k}}\}_{2\le k\le s-1}\}\in\ba_{s-2}^1 : i_0, i_1\in \tilde I, \ \ i_0<i_1<i_2\},
\\
\ba_s^5=\{\{\gamma_{1,n-1},\gamma_{1,n},\delta_{i_1},\{\beta _{i_{k},j_{k}}\}_{2\le k\le s-2}\}: \{\beta _{i_{k},j_{k}}\}_{2\le k\le s-2}\}\in\ba_{s-3}^1 :  i_1\in \tilde I, \ \ <i_1<i_2\},\\
\ba_s^6=
 \{\{\tilde{\gamma}_{i_0}, \{\beta _{i_{k},j_{k}}\}_{1\le k\le s-1}\}: \{\beta _{i_{k},j_{k}}\}_{1\le k\le s-1}\}\in\ba_{s-1}^1,\ i_0\in\tilde I, i_0<i_1
 \},
\end{gather*}
where $\tilde{\gamma}_{i_0}\in\{\gamma_{i_0,n-1},\gamma_{i_0,n},\delta_{i_0}\}$ and $\gamma_{i_0}\in\{\gamma_{i_0,n-1},\gamma_{i_0,n}\}$.
This gives, $$\#\ba_s = \left(
\binom{n-2}{2s}+ (\binom{n-3}{2s-2}+\binom{n-3}{2s-1})+3\binom{n-2}{2s-1}+\binom{n-3}{2s-4}+
2\binom{n-2}{2s-1}+\binom{n-3}{2s-5}\right),$$ and hence
$\sum_s \#\ba_s=2^n$.\hfill\qedsymbol

\end{prop}
In the corollary we understand that all the missing cases correspond to the empty set.
\begin{cor}
Suppose that $J$  is a subset of $I$. Then $\ba_{s,J}=\sqcup_{p=1}^6\ba^p_{s,J}$,
where  $\ba_{s,J}^p=\ba_{s,J}^{p,1}\sqcup\ba_{s,J}^{p,2}$ are given by,
\begin{gather*}\ba_{s,J}^{1,1} =\{ \{ \beta _{i_{k},j_{k}}\} _{1\leq k\leq
s}\in \ba_{s}^{1}:i_{k},j_{k}\notin J\}  \\
\ba_{s,J}^{1,2} =\{ \{ \beta _{i_{k},j_{k}}\} _{1\leq k\leq
s}\in \ba_{s}^{1}:\{ \beta _{i_{k},j_{k}}\} _{1\leq k\leq s-1}\in
\ba_{s-1,J}^{1,1},\text{ }i_{s}\in J,\text{ }j_{s}=i_{s}+1\notin J,\} \\ \ba_{s,J}^{2,q} =
 \{ \alpha _{1,l},\{ \beta _{i_{k},j_{k}}\} _{1\leq
k\leq s-1}\}\in\ba_{s}^2 :\{ \beta _{i_{k},j_{k}}\} _{1\leq k\leq
s-1}\in \ba_{s-1,J}^{1,q},\text{ } \ell\notin J\}\  \rm{if}\ \ 1\notin J,\\
\ba_{s,J}^{3,q}=\{\{\gamma_{1,n-1},\gamma_{1,n},\{\beta _{i_{k},j_{k}}\}_{1\le k\le s-2}\}\in\ba_{s}^3: \{\beta _{i_{k},j_{k}}\}_{1\le k\le s-2}\}\in\ba_{s-2}^{1,q}\ \}\ \ {\rm{if}}\ 1, n-1, n\notin J,\\
\ba_{s,J}^{4,q} =\{
 \{ \gamma _{i_0},\delta _{i_1},\{ \beta
_{i_{k},j_{k}}\} _{1\leq k\leq s-2}\}\in\ba_s^4 :\{ \beta
_{i_{k},j_{k}}\} _{1\leq k\leq s-2}\in \ba_{s-2,J}^{1,q},\text{ }%
i_0,i_1\notin J\} , \ \  n-1,n\notin J,\\%
  \ba_{s,J}^{5,q} =\{
 \{ \gamma _{1,n-1},\gamma _{1,n},\delta _{i_1},\{ \beta
_{i_{k},j_{k}}\} _{1\leq k\leq s-3}\}\in\ba_s^5 :\{ \beta
_{i_{k},j_{k}}\} _{1\leq k\leq s-3}\in \ba_{s-3,J}^{1,q},\text{ }i_1\notin
J\},\  \ 1,n-1,n\notin J
\end{gather*} where $q\in\{1,2\}$ and $\gamma _{i_0}\in \{ \gamma _{l,n-1},\gamma_
{l,n}\}$
and $\ba_{s,J}^{6,q}$ is given by
\begin{gather*}\ba_{s,J}^{6,q}=
 \{\{\tilde{\gamma}_{i_0}, \{\beta _{i_{k},j_{k}}\}_{1\le k\le s-1}\}\in\ba_s^6: \{\beta _{i_{k},j_{k}}\}_{1\le k\le s-1}\}\in\ba_{s-1}^{1,q}\ i_0\notin J\},
\end{gather*}provided that  $n-1\notin J$ if $\tilde{\gamma}_{i_0}=\gamma_{i_0,n-1}$ and $n\not\in J$ if $\tilde{\gamma}_{i_0}=\gamma_{i_0}$ and $n-1,n\notin J$ if $\tilde{\gamma}_{i_0}=\delta_{i_0}$.
\hfill\qedsymbol

\end{cor}

\section{Irreducible ad-nilpotent ideals}
We use the notation of section one freely.  Thus, $\lie g$ is a simple Lie algebra, $\lie h$ is a fixed Cartan subalgebra of $\lie g$ and $R$ is the corresponding set of roots. However, we do not fix a set of simple roots or a positive system once and for all, rather we choose it depending on our situation. Anytime that we do pick a particular set of simple roots, then we understand that the associated data $Q^+$, $P^+$ etc. are all as defined in section one.
 \subsection{} Let $\lie p$ be a parabolic subalgebra of $\lie g $ containing $\lie h$ and let $\lie p=\lie m\oplus \lie u$ be its Levi decomposition, where $\lie m$ is the  reductive part  and $\lie u$ is the unipotent radical of $\lie p$. The restriction of the adjoint action of $\lie g$  to $\lie m$ induces on $\lie u$ the   structure  of an $\lie m$--module.  The centre  $\lie z$ of $\lie m$  acts semisimply on $\lie u$. An unpublished result of Kostant (a proof can be found in \cite{Jos}, \cite{Wolf}) is that the distinct  $\lie z$--eigenspaces of  $\lie u$ are irreducible $\lie m$--modules. These can be described as follows.
 Choose a positive root system $R^+$ and a proper subset $J$ of the simple roots so that $$\lie p=\lie h\oplus_{\alpha\in R^+}\lie g_\alpha\oplus_{\alpha\in R^+(J)}\lie g_{-\alpha}, \ \lie m= \lie h\oplus_{\alpha\in R^+(J)}\lie g_{\pm\alpha},\ \ \lie u=\oplus_{\{\alpha\in R^+\setminus R^+(J)\}}\lie g_\alpha.$$ Let $\sim$ be the equivalence relation on $R^+$ given by $\alpha\sim\beta$ iff $d_j(\alpha-\beta)\ne 0$ implies $j\in J$. An irreducible $\lie m$--submodule of $\lie u$  is just the direct sum of root spaces which lie in a fixed equivalence class.   The following is an easy exercise, but we include a proof for completeness.
 
\begin{lem}\label{unique} There exists a  unique  ad--nilpotent ideal $\lie i_0$  of $\lie p$  which is irreducible as an $\lie m$--module. In fact, $\lie i_0=\oplus_{\alpha\sim\theta}\lie g_\alpha$ and in particular  $\lie i_0$ is abelian.

\end{lem}

\begin{pf} Assume that we have fixed a system as in the discussion preceding the statement of the Lemma.  Let  $\alpha, \gamma\in R^+$ with $\alpha+\gamma\in R^+$, then $\alpha + \gamma\le\theta$.  If  $\alpha\sim\theta$ then $d_j(\theta-\alpha)\ne 0$ only if $j\in J$ and hence $d_j(\gamma)\ne 0$ only if $j\in J$. This proves that
 that $\gamma\in R^+(J)$ and hence $\alpha+\gamma\sim\theta$. It also proves that if $\alpha\sim\theta$ and $\gamma\sim\theta$ then $\alpha+\gamma\notin R^+$. Since $\alpha-\beta\sim\theta$ for all $\alpha \sim \theta$ and $\beta\in R^+(J)$, it follows that $\lie i_0$ is an abelian ideal of $\lie p$. The fact that it is irreducible as a $\lie m$--module is a special case of the result of Kostant mentioned above. If $\lie i$ is any other ad--nilpotent ideal in $\lie p$, then $x_\theta\in\lie i$ and hence $\lie i$ contains the irreducible $\lie m$--module generated by $x_\theta$, i.e $\lie i_0\subset\lie i$. \hfill\qedhere

\end{pf}

 We call $\lie i_0$ the irreducible ad--nilpotent ideal of $\lie p$. In other words, given any parabolic subalgebra $\lie p$  of $\lie g$ containing $\lie h$ we can associate to it, canonically a certain subset of roots of $R$: namely the roots that determine the irreducible ad--nilpotent ideal of $\lie p$. It is not difficult to see that different parabolic subalgebras could give rise to the same set of roots.
 In this section we give a necessary and sufficient condition on a  subset $S$  of $R$ to determine the unique irreducible ideal $\lie i_S$ in some parabolic subalgebra of $\lie g$ which contains $\lie h$.
 As a corollary we prove that given such a subset there exists a parabolic $\lie p_S$  whose unique irreducible ad--nilpotent ideal is given by $S$ and  if $\lie p$ is another parabolic whose irreducible ad nilpotent ideal is given by $S$ then $\lie p\subset \lie p_S$.
 As a further consequence of our result we write down explicitly all subsets of $R$  for all the classical Lie algebras which satisfy this condition. As explained in the introduction, our motivation for describing such sets comes from the results of \cite{CG}, where to each such set,  we construct a finite and infinite-dimensional Koszul algebra of global dimension equal to the cardinality of the set and we expect that this description will be useful in the further study of those algebras.

\subsection{} Given a subset $S$ of $R$, let  $\bz_+S\subset Q$ be the $\bz_+$--span of elements of $S$ and define $\rho_S\in P$ by
  $2\rho_S=\sum_{\alpha\in S}\alpha.$ Given $\lambda\in P$ set $$\max\lambda=\max\{(\lambda,\alpha):\alpha\in R\},\ \qquad S(\lambda)=\{\alpha\in R:(\lambda,\alpha)=\max\lambda\}.$$
 Define a function $\bod: Q^+\to\bz_+$ by, $\bod(\eta)=\ \min\{\sum_{\alpha\in R}m_\alpha: \eta=\sum_{\alpha\in R} m_\alpha\alpha\}.$
Our main result is the following.
\begin{thm}\label{main2} Let $S\subset R$. The following conditions are equivalent.
\begin{enumerit} \item  $S=S(\rho_S),\ \ \max\rho_S>0$
\item $S=S(\lambda)$ for some $\lambda\in P$ with $\max\lambda>0$.
\item If $\eta\in\bz_+S$ and  $\eta=\sum_{\alpha\in R^+}m_\alpha\alpha$ for some $m_\alpha\in\bz_+$ then $\bod(\eta)=\sum_{\alpha\in R} m_\alpha$ iff $ m_\alpha=0$ for all $\alpha\notin S.$
\item For all
$\alpha,\beta\in S$ we have  $\alpha+\beta\notin R$ and if $\gamma\in R\setminus S$, then $(\gamma+R)\cap (S+S)=\emptyset$.

\end{enumerit}
\end{thm}
\subsection{} Before proving the theorem, we establish the following corollary. Given a non--zero  $\lambda\in P$  let $\lie p_\lambda$ be the corresponding parabolic subalgebra of $\lie g$, i.e, $$\lie p_\lambda=\lie m_\lambda\oplus\lie u_\lambda,\qquad \lie m_\lambda=\lie h\bigoplus_{\{\alpha: (\lambda,\alpha)=0\}}\!\!\!\!\!\!\!\lie g_{\pm\alpha},\qquad \ \lie u_\lambda=\bigoplus_{\scriptsize{\{\alpha: (\lambda,\alpha)>0)\}}}\!\!\!\!\!\!\!\lie g_\alpha.$$

\begin{cor} Let $S$ satisfy one of the equivalent conditions of the theorem. Then,
\begin{enumerit}
\item[(i)]
The unique irreducible ad--nilpotent  ideal of $\lie p_{\rho_S}$ is given by $S$.
\item[(ii)]
If $\lie p$ is another parabolic subalgebra of $\lie g$ containing $\lie h$ whose  unique irreducible ad--nilpotent  ideal is given by $S$, then $\lie p\subset\lie p_{\rho_S}$.\end{enumerit}
\end{cor}
\begin{pf} Since $S=S(\rho_S)$, to prove (i) it suffices by Lemma \ref{unique} to prove that $$\lie i_{\rho_S}= \oplus_{\alpha\in R^+:(\rho_S,\alpha)=\max\rho_S}\lie g_\alpha$$ is an ideal which is irreducible as an $\lie m_{\rho_S}$--module. The fact that $\lie i_{\rho_S}$
is an ideal in $\lie p_S$ is trivial. Choose a positive root system $R^+$  so that $\rho_S$ is dominant integral with respect to $R^+$, in which case $S\subset R^+$ and $\theta\in S$. If  $\alpha\in S$ then $(\rho_S,\theta-\alpha)=0$ and since $\theta-\alpha=\sum_{k=1}^p\beta_k$ for some $\beta_1,\cdots,\beta_p\in R^+$ it follows that $(\rho_S,\beta_k)=0$ for all $1\le k\le p$, which proves that $\alpha\sim\theta$. Hence $\lie i_{\rho_S}$ is contained in the unique irreducible ad--nilpotent ideal of $\lie p_{\rho_S}$ and part (i) is proved.

To prove (ii) suppose that $\lie p$ is a parabolic subalgebra and that $\lie i_{\rho_S}\subset \lie p$ is the unique irreducible ad--nilpotent ideal in $\lie p$.  It suffices to show that $(\rho_S,\alpha)\ge 0$ for all $\alpha\in R$  with $\lie g_\alpha\subset\lie p$ and that $(\rho_S,\alpha)=0$ if $\lie g_{\pm\alpha}\subset \lie p$.  Suppose that for some $\beta\in S$ and $\alpha\in R$ we have $(\beta,\alpha)<0$. Then  $\beta+\alpha\in R$  and if $\lie g_\alpha\subset \lie p$ then  $\lie g_{\alpha+\beta}\subset\lie i_{\rho_S}$ i.e, $\alpha+\beta\in S$.
Hence $$(\rho_S,\beta+\alpha)=(\rho_S,\beta)=\max\rho_S,$$ which implies $(\rho_S,\alpha)=0$.  On the other hand if $\alpha\in R$ is such that $(\beta,\alpha)\ge 0$ for all $\beta\in S$, then clearly $(\rho_S,\alpha)\ge 0$. It remains to prove that $(\rho_S,\alpha)=0$ if $\lie g_{\pm\alpha}\subset \lie p$. But this too is clear since either $(\beta,\alpha)=0$ for all $\beta\in S$, or there exists $\beta\in S$ such that either  $(\beta,\alpha)<0$ or $(\beta,-\alpha)<0$ and in either case we have seen that $(\rho_S,\alpha)=0$.

\end{pf}

\subsection{} The rest of the section is  devoted to to proving the theorem. It is trivially true  that (i) implies (ii).
 To prove that (ii) implies (iii)   suppose that $\eta\in \bz_+S(\lambda)$ and suppose that we have  $$\eta=\sum_{\alpha\in R}k_\alpha\alpha=\sum_{\beta\in S(\lambda)}m_\beta\beta,\ \ k_\alpha, m_\alpha\in\bz_+,$$ with $\bod(\eta)=\sum_{\alpha\in R}k_\alpha\le \sum_{\beta\in S}m_\beta$. Then, $$(\lambda,\eta)=\max\lambda\sum_{\beta\in S}m_\beta=\sum_{\alpha\in  R}k_\alpha(\lambda,\alpha)\le \max\lambda \sum_{\alpha\in R}k_\alpha=\max\lambda\ \bod(\eta).$$ Since $\max\lambda>0$, it follows that $\sum_{\beta\in S}m_\beta=\bod(\eta)$ and also that $\sum_{\alpha\in R}k_\alpha(\max\lambda-(\lambda,\alpha))=0$. Since $k_\alpha\in\bz_+$ and $\max\lambda-(\lambda,\alpha)\ge 0$ it follows that $k_\alpha=0$ if $\alpha\not\in S$ and we are done.

Assume that $S$ satisfies the conditions of (iii). Suppose that $\alpha,\beta\in S$ is such that $\alpha+\beta\in R$. Then we have $\bod(\alpha+\beta)=1$ and $\bod(\alpha+\beta)=2$ which is absurd. Next suppose that $\gamma\notin S$ and that there exists $\delta\in R$ such that $\gamma+\delta=\alpha+\beta$ for some $\alpha,\beta\in S$. Since $\alpha+\beta\in \bz_+S$ we have  $\bod(\alpha+\beta)=2=\bod(\gamma+\beta)$ and so  $\gamma,\delta\in S$ by the hypothesis of (iii).

\subsection{}  We shall need the following result  to prove that
(iv) implies (i). Note that the condition that $\alpha+\beta\notin R$ if $\alpha,\beta\in S$ implies that $(\alpha,\beta)\ge 0$ for all $\alpha,\beta\in S$. We shall use this remark freely throughout the rest of the section.

\begin{lem}\label{lemone}Let $\alpha,\beta\in R$. Assume that $\beta$ is a long root and that $(\alpha,\beta)=0$.
  There exists
 $\gamma,\gamma'\in R$ with $\gamma\notin\{\alpha,\beta\}$  such that $\alpha+\beta=\gamma+\gamma'$.
\end{lem}
\begin{pf}
Let $\{\alpha_i:i\in I\}$ be a set of simple roots and assume without loss of generality that $\alpha$ is
simple, say $\alpha=\alpha_{i_0}$ and set $I_0=\{i_0\}$.
 For $k\ge 1$, define subsets $I_k$ recursively by,
$$\ \ I_k=\{i\in I: (\alpha_i,\alpha_j)<0\ \ {\text{for some}} \ j\in I_{k-1}\}. $$  Clearly  $I=\cup_{k\ge 0} I_k,$ and we pick $m$ minimal so that there exists $i_m\in I$ with $(\beta,\alpha_{i_m})\ne 0$.
Since  $i_m\notin I_{m-1}$ we pick
$i_{m-1}\in I_{m-1}$ with $(\alpha_{i_m}, \alpha_{i_{m-1}})<0$.
Again, $i_{m-1}\notin I_{m-2}$ since otherwise we would have
$i_m\in I_{m-1}$. In other words, we can choose elements
 $i_k\in I_k\setminus I_{k-1}$, with $1\le k<m$ such that $(\beta, \alpha_{i_k})=0$  and$(\alpha_{i_k}, \alpha_{i_{k+1}})<0$.
Hence $\alpha_{i_p}+\cdots +\alpha_{i_k}\in R$ for all $0\le p\le k\le m$.  If $(\beta,\alpha_{i_m})<0$, then
 $(\beta,\alpha_{i_p}+\cdots +\alpha_{i_m} )<0$ for all $0\le p\le m$
and so $\gamma= \alpha_{i_0}+\alpha_{i_1}+\cdots+\alpha_{i_m}+\beta\in R$. Since $\gamma-\alpha_{i_0}, \gamma-\beta\in R$, they are in particular non--zero and the Lemma follows, from
$$\alpha_{i_0}+\beta= (\alpha_{i_0}+\alpha_{i_1}+\cdots+\alpha_{i_m}+\beta)- (\alpha_{i_1}+\cdots+\alpha_{i_m}).$$

If $(\beta,\alpha_{i_m})>0$, then
$\gamma=\alpha_{i_0}+\alpha_{i_1} +\cdots +\alpha_{i_m}\in R$ and also, $\gamma'=\beta- (\alpha_{i_1} +\cdots +\alpha_{i_m})\in R. $
 It remains to prove that
 $\gamma\notin\{\alpha,\beta\}$. It is clear that
 $\gamma\ne\alpha_{i_0}$ since $\gamma-\alpha_{i_0}\in R$. Suppose
 that $\gamma=\beta$. Then we have
 $$(\beta,\beta)=(\beta,\alpha_{i_0}+\alpha_{i_1} +\cdots +\alpha_{i_m})=(\beta,\alpha_{i_m}).$$ Since $\beta$ is a long root this  implies
 that  we must have  $\beta=\alpha_{i_m}=\gamma$. But this is impossible since $\alpha_{i_0}+\alpha_{i_1} +\cdots +\alpha_{i_{m-1}}\in R$ and the proof of the Lemma is complete.

\end{pf}

\subsection{} We first prove that (iv) implies (i) when  $S$ is a subset of the  of long roots in $R$. For $\alpha\in S$, set $$\alpha^\perp=\{\gamma\in S: (\alpha,\gamma)=0\}.$$
We claim that for all $\alpha,\beta\in S$, we have,
\begin{equation}\label{one} \#\alpha^\perp=\#\beta^\perp.\end{equation}
Assuming the claim, the proof that $S=S(\rho_S)$ is completed as follows. Since $S$ consists of long roots, we  assume that $(\delta,\delta)=2$ if $\delta\in S$ and also that $(\delta,\beta)\le 1$ for all $\beta\in R\setminus\{\delta\}$. If $\gamma\in S$ and $\gamma\ne \delta$, then we have  $0\le (\gamma,\delta)\le 1$ and so
 $$(\rho_S,\gamma)=\sum_{\delta\in S}(\delta,\gamma)= \#S+1-\#\gamma^\perp,\ \ \gamma\in S.$$  Suppose now that $\gamma'\notin S$. If $(\gamma,\delta)\le 0$
for all $\delta\in S$, then $(\rho_S,\gamma)\le 0$ and hence $\gamma'\notin S(\rho_S)$.  Suppose that  $(\gamma',\delta)>0$
 for some $\delta\in S$. If  $\delta_1\in \delta^\perp$ is such that
$(\delta_1,\gamma')>0$ then we would have that $\delta_1-\gamma'+\delta\in R$. Since $$\delta_1+\delta=(\delta_1-\gamma'+\delta)+\gamma$$ it would follow that  $\gamma'\in S $ which is a contradiction. Hence, $(\gamma',\delta^\perp)=0$. But now, we have
$$(\rho_S,\gamma')=\sum_{\delta_1\in \delta^\perp}(\delta_1,\gamma')+
\sum_{\delta_1\in S\setminus \delta^\perp)}(\delta_1,\gamma')\le
 \#S-\#\delta^\perp<\#S+1-\#\delta^\perp=\rho_S(\delta).$$ Hence, $\gamma'\notin S(\rho_S)$ proving that $S=S(\rho_S)$.

It suffices to prove \eqref{one} when  $(\alpha,\beta)>0$; for,
  if $(\alpha,\beta)=0$, then by Lemma \ref{lemone}, we can choose $\gamma,\gamma'\in R$ such that $\alpha+\beta=\gamma+\gamma'$.
  By the conditions on $S$, this means that $\gamma,\gamma'\in S$.
  Using the fact that $S$ consists of long roots, we see also that $$ (\alpha,\gamma)+(\beta,\gamma)=(\gamma,\gamma)+(\gamma,\gamma')\ge 2,$$
  which in turn implies that $(\alpha,\gamma)>0$ and $(\beta,\gamma)>0$, which gives $$\#\alpha^\perp=\#\gamma^\perp=\#\beta^\perp.$$
Assume now that $(\alpha,\beta)>0$. Let $s_\alpha,s_\beta$ be the reflections in $W$ corresponding to the roots $\alpha$ and $\beta$. It suffices to prove that
\begin{equation}\label{two} s_\alpha s_\beta(\alpha^\perp)\subset \beta^\perp.\end{equation}
If $\gamma\in\alpha^\perp$, it is easy to calculate that
 $$s_\alpha s_\beta(\gamma)=\begin{cases}\gamma,\ \
  (\gamma,\beta)=0,\\ \gamma-\beta+\alpha,\ \ \ (\gamma,\beta)=1.\end{cases}.$$ In the first case, $\gamma\in S(\beta)$.
   In the second case we have $\alpha+\gamma=(\gamma-\beta+\alpha)+\beta$ which implies that $\gamma-\beta+\alpha\in S$.
   Since $S$ consists of long roots, we get $(\gamma-\beta+\alpha,\beta)=0$ and  \eqref{two} follows.

\subsection{} It remains to consider the case when  $S$ contains  a short root. It is clear that if $S$ satisfies the conditions of Theorem 2 (iv), then so does $wS$ for all $w\in W$. Hence, we fix a set of simple roots and  assume that either $S$ contains a short simple root or the highest short root.
Suppose first that $\lie g$ is of type $G_2$, and let $\alpha_1,\alpha_2$ be the simple roots with $\alpha_2$ being short. Assume that $\alpha_2\in S$. Since $2\alpha_2=(\alpha_1+2\alpha_2)-\alpha_1$ and $\alpha_1+2\alpha_2\in R$,  it follows that $\alpha_1+2\alpha_2\in S$, but this is impossible since $\alpha_1+3\alpha_2\in R$. Hence, the result is vacuously true for $G_2$.
Suppose now that $\lie g$ is of type $F_4$ and suppose that $S$ contains the highest short root. Then, a similar argument proves
$$ S=\{ \alpha\in R^+: d_4(\alpha)=2\},\ \ \rho_S = 7\omega_4.$$

The following  proposition, which lists all sets $S$ satisfying the hypothesis of Theorem 2(iv) for the classical Lie algebras, also completes the proof in the case when $S$ contains a short root. For the purposes of this proposition, it is convenient to use a slightly different notation for the roots. Let $\{\epsilon_i$\ : \ $1\le i\le n+1$\} be a basis for $\br^{n+1}$.  For $\lie g$ of type $A_n$, we set $\alpha_i = \epsilon_i - \epsilon_{i+1}$, $1\le i\le n$. For $\lie g$ of types $B_n,\ C_n$, and $D_n$, we set $\alpha_i = \epsilon_i - \epsilon_{i+1}$ if $1 \le i \le n-1$. Finally, we set $\alpha_n = \epsilon_n,\ \alpha_n = 2\epsilon_n$, and $\alpha_n = \epsilon_{n-1} + \epsilon_n$, for $\lie g$ of type $B_n,\ C_n$, and $D_n$, respectively.

\begin{prop} Let $S\subset R$ be such that, if $\alpha,\beta\in S$, then $\alpha+\beta\notin R$ and, if $\gamma\in R\setminus S$, then $(\gamma+R)\cap (S+S)=\emptyset$.
\begin{enumerit}
\item If $\lie g$ is of type $A_n$, there exist  disjoint subsets~$\boi$, $\boj$ of $\{1,\cdots ,n+1\}$ such that
$$S=\{\epsilon_{i}-\epsilon_j: i\in\boi, j\in\boj\},\ \ \qquad \rho_S=|\boj|\sum_{p=1}^{|\boi|}\epsilon_{i_p}-|\boi|\sum_{p=1}^{|\boj|}\epsilon_{j_p}.$$

    \item If $\lie g$ is of type $C_n$, there exist disjoint subsets $\boi$ and $\boj$ of $\{1,\cdots,n\}$ such that
     \begin{gather*}S=\{\epsilon_{i_1}+\epsilon_{i_2}: i_1,i_2\in\boi\}\cup\{-(\epsilon_{j_1}+\epsilon_{j_2}): \ j_1,j_2\in\boj\}\cup\{\epsilon_i-\epsilon_j: i\in\boi, j\in\boj\},\\ \rho_S =(|\boi|+|\boj|+1)(\sum_{i\in\boi}\epsilon_i-\sum_{j\in\boj}\epsilon_j).\end{gather*}
    \item If $\lie g$ is of type $B_n$ then either
    \begin{enumerit} \item[(a)] there exist  disjoint subsets $\boi$ and $\boj$ of $\{1,\cdots,n\}$ such that \begin{gather*}S=\{\epsilon_i-\epsilon_j: i\in\boi, j\in\boj\}\ \cup\ \{\epsilon_k+\epsilon_\ell: k,\ell\in\boi, k\ne \ell \}\cup\{-(\epsilon_k+\epsilon_\ell):\ k, \ell\in\boj,\ \ k\ne \ell\},\\ \rho_S=(|\boi|+|\boj|-1)(\sum_{i\in\boi}\epsilon_i-\sum_{j\in\boj}\epsilon_j),\end{gather*}
    \item[(b)] or there exists $i\in I$ such that \begin{gather*}  \pm S = \{\epsilon_i\}\cup\{\epsilon_i\pm\epsilon_j\ :\ i\ne j, 1\le j\le n\},\\ \rho_S=\pm(2n-1)\epsilon_i.\end{gather*}
    \end{enumerit}
\item If $\lie g$ is of type $D_n$ then either
    \begin{enumerit} \item[(a)] there exist  disjoint subsets $\boi$ and $\boj$ of $\{1,\cdots,n\}$ such that \begin{gather*}S=\{\epsilon_i-\epsilon_j: i\in\boi, j\in\boj\}\ \cup\ \{\epsilon_k+\epsilon_\ell: k,\ell\in\boi, k\ne \ell \}\cup\{-(\epsilon_k+\epsilon_\ell):\ k, \ell\in\boj,\ \ k\ne \ell\},\\\rho_S=(|\boi|+|\boj|-1)(\sum_{i\in\boi}\epsilon_i-\sum_{j\in\boj}\epsilon_j),\end{gather*}
    \item[(b)] or there exists $i\in I$ such that $$ \pm S=\{\epsilon_i\pm\epsilon_j:i\ne j, 1\le j\le n\},\ \ \rho_S= \pm(2n-2)\epsilon_i.$$
    \end{enumerit}

    \end{enumerit}
    \end{prop}
    \begin{pf} The proof is elementary.  For instance, suppose that $\lie g$ is of type $A_n$ and that  $\alpha= \epsilon_i-\epsilon_j$ and $\beta=\epsilon_k-\epsilon_\ell$ are in $S$ and $\alpha\ne\beta$.  Then $(\alpha,\beta)\ge 0$ implies  $\ell\ne i$ and $j\ne k$. Since
   $\alpha+\beta= \epsilon_i-\epsilon_\ell+ \epsilon_k-\epsilon_j$ we must have $\epsilon_i-\epsilon_\ell, \epsilon_k-\epsilon_j\in S$.  This proves (i). We omit the equally simple proofs  of the other cases.
    \end{pf}


\begin{thebibliography}{CPS2}


\bibitem{B} N. Bourbaki, {\em Groupes et alg\`ebres de {L}ie,
Chapitres 4, 5 et 6,} Hermann, Paris 1968, Masson, Paris 1981;
{\em {L}ie groups and {L}ie algebras, Chapter 4--6,}
Translated from the 1968 French original by Andrew Pressley,
Springer, Berlin 2002



\bibitem{CP1} P. Cellini and P. Papi, {\em $ad$-nilpotent ideals of a Borel subalgebra}, J. Algebra {\bfseries 225} (2000), 130--141

\bibitem{CP2} \bysame, {\em $ad$-nilpotent ideals of a Borel subalgebra II}, J. Algebra {\bfseries 258} (2002), 112--121

\bibitem{CP3} \bysame, {\em $ad$-nilpotent ideals of a Borel subalgebra III}, 2003, {\tt math.RT/0303065v2}

\bibitem{CP4} \bysame, {\em The structure of total reflection orders in affine root systems}, J. Algebra {\bfseries 205} (1998), 207--226

\bibitem{CG} V. Chari and J. Greenstein, {\em A Family of Koszul algebras arising from finite-dimensional representations of simple Lie algebras}, 2008, {\tt math.RT/0808.1463v1}

\bibitem{Hum} J. E. Humphreys, {\em Introduction to Lie algebras and representation theory}, Springer, New York, 1978

\bibitem{Jos} A. Joseph, {\em Orbital varieties of the minimal orbit},Ann. Ec. Norm. Sup. {\bfseries 31} (1998),
17--45



\bibitem{Ko} B. Kostant, {\em The Set of
Abelian ideals of a Borel Subalgebra, Cartan Decompositions, and Discrete Series
Representations}, Internat. Math. Res. Notices {\bfseries 5} (1998), 225--252

\bibitem{KOP} C. Krattenthaler, L. Orsina and P. Papi, {\em Enumeration of
$ad$-nilpotent $\frak b$-ideals for simple Lie algebras}, Adv. Appl. Math. {\bfseries 28} (2002), 478--522, {\tt math.RA/0011023}

\bibitem{OP} L. Orsina and P. Papi, {\em Enumeration of $ad$-nilpotent ideals of a Borel subalgebra in type A by class of nilpotence}, C.R. Acad. Sci. Paris S\`er. I Math. {\bfseries 330} (2000), no.8, 651--655

\bibitem{Pan} D. Panyushev, {\em Abelian ideals of a Borel subalgebra and long positive roots}, IMRN {\bfseries 35} (2003), 1889--1913

\bibitem{Pan2} \bysame, {\em 	$ad$-nilpotent ideals: generators and relations}, J. Algebra {\bfseries 274} (2004), 822--846


\bibitem{Rig} C. Righi, {\em Ad-Nilpotent Ideals of a Parabolic Subalgebra}, J. Algebra {\bfseries 319} (2008), 1555--1584, {\tt math.RT/0701679v2}

\bibitem{Rig2} \bysame, {\em Number of '`udu'` of a Dyck Path and $ad$-nilpotent Ideals of Parabolic Subalgebras of $sl_{l+1}(\cal{C})$}, 2008, {\tt math.RA/0803.0267v1}

\bibitem{Sut} R. Suter, {\em Abelian ideals in a Borel subalgebra of a complex simple Lie algebra}, Inventiones mathematicae {\bfseries 156} (2004), 175--221, {\tt math.RT/0210463v1}

\bibitem{Wolf} J. Wolf, {\em Sapces of Constant Curvature}, McGraw-Hill, New York, 1962


\end{thebibliography}
\end{document}